\documentclass[leqno,12pt]{amsart}

\usepackage{latexsym,esint}
\usepackage{color}

\theoremstyle{plain}

\def\endproof{\hspace*{\fill}\mbox{\ \rule{.1in}{.1in}}\medskip }

\newtheorem{theorem}{Theorem}[section]

\newtheorem{lemma}[theorem]{Lemma}

\theoremstyle{definition}

\numberwithin{equation}{section}

\begin{document}
\title[viscoelasticity with physical viscosity]
{A local existence result for system \\ of viscoelasticity with physical viscosity}
\author{Marta Lewicka}
\address{Marta Lewicka,  University of Pittsburgh, Department of Mathematics, 
301 Thackeray Hall, Pittsburgh, PA 15260, USA }
\email{lewicka@math.pitt.edu}
\author{Piotr B. Mucha}
\address{Piotr B. Mucha, Institute of Applied Mathematics and Mechanics, 
 University of Warsaw, ul. Banacha 2, 02097 Warszawa, Poland}
\email{p.mucha@mimuw.edu.pl}

\begin{abstract}
We prove the local in time existence of
the classical solutions to the system of equations of isothermal
viscoelasticity with clamped boundary conditions. 
We deal with a general form of viscous stress tensor
$\mathcal{Z}(F,\dot F)$, assuming a Korn-type condition on its 
derivative $D_{\dot F}\mathcal{Z}(F, \dot F)$. This condition
is compatible with the balance of angular momentum, frame invariance
and the Claussius-Duhem inequality. We give examples of linear and
nonlinear (in $\dot F$) tensors $\mathcal{Z}$ satisfying these
required conditions.
\end{abstract}

\maketitle

\section{Introduction and the main results}

In this paper, we are concerned with the local in time existence of
the classical solutions to the system of equations of isothermal
viscoelasticity. The system we study is given through the
balance of linear momentum:
\begin{equation}\label{v_eq}
\xi_{tt} - \mbox{div}\Big(DW(\nabla\xi) + \mathcal{Z}(\nabla\xi,
\nabla\xi_t)\Big)=0 \quad \mbox{ in } \Omega\times \mathbb{R_+},
\end{equation}
and it is subject to initial data:
\begin{equation}\label{initial}
 \xi(0,\cdot)=\xi_0 ~~~\mbox{ and } ~~~\xi_t(0,\cdot)=\xi_1 \quad \mbox{ in } 
\Omega,
\end{equation}
the clamped boundary conditions:
\begin{equation}\label{bc}
 \xi(\cdot, X)=X \qquad \forall X\in\partial \Omega, 
\end{equation}
and the non-interpenetration ansatz:
\begin{equation}\label{nonin}
\det \nabla \xi > 0 \quad \mbox{ in } \Omega.
\end{equation}
Here, $\xi:\Omega\times\mathbb{R_+}\longrightarrow \mathbb{R}^n$
denotes the deformation of a reference configuration
$\Omega\subset\mathbb{R}^n$ which models
a viscoelastic body with constant temperature and density.
A typical point in $\Omega$ is denoted by $X$, and 
the deformation gradient, the velocity and velocity gradient are given as:
$$F = \nabla\xi\in\mathbb{R}^{n\times n}, \quad
v=\xi_t\in\mathbb{R}^n, \quad Q=\nabla\xi_t = \nabla v = F_t\in\mathbb{R}^{n\times n}.$$
In (\ref{v_eq}) the operator ${div}$ stands for the spacial divergence of an appropriate field.
We use the convention that the divergence of a matrix field is taken
row-wise. In what follows, we shall also use the matrix norm
$|F|=(\mbox{tr}(F^TF))^{1/2}$, which
is induced by the inner product: $F_1:F_2 = \mbox{tr}(F_1^TF_2)$.
To avoid notational confusion, we will often write $\langle
F_1:F_2\rangle$ instead of $F_1:F_2$.

\subsection{The elastic energy density $W$.}
The mapping $DW:\mathbb{R}^{n\times n}\longrightarrow \mathbb{R}^{n\times n}$
in (\ref{v_eq}) is the Piola-Kirchhoff stress tensor which, in agreement with
the second law of thermodynamics \cite{Dafbook},
is expressed as the derivative of an elastic energy density
$W:\mathbb{R}^{n\times n}\longrightarrow \overline{\mathbb{R}}_+$.

The principles of material frame invariance, material consistency, and normalisation 
impose the following conditions on $W$, valid for all
$F\in\mathbb{R}^{n\times n}$ and all proper rotations $R\in SO(n)$:
\begin{equation}\label{elastic_dens}
\begin{minipage}{14cm}
\begin{itemize}
\item[(i)] $W(RF) = W(F),$
\item[(ii)] $W(F)\to +\infty \quad \mbox{ as } \det F\to 0$,
\item[(iii)] $W(\mbox{Id}) = 0$.
\end{itemize}
\end{minipage}
\end{equation}
Examples of  $W$ satisfying the above conditions are: 
\begin{equation*}
\begin{split} 
W_1(F) & = |(F^TF)^{1/2} - \mbox{Id}|^2 + |\log \det F|^q, \\
W_2(F) & =  |(F^TF)^{1/2} - \mbox{Id}|^2 + \left|\frac{1}{\det F} - 1\right|^q 
\mbox{ for } \det F>0,
\end{split}
\end{equation*}
where $q>1$ and $W$ is intended to be $+\infty$ if $\det F\leq 0$ \cite{MS}.
Another case-study  example, satisfying (i) and (iii) is:
$W_0(F)=|F^TF-\mbox{Id}|^2 = (|F^T F|^2 - 2|F|^2+n)$.

We will assume that $W$ is smooth in a neighborhood of $SO(n)$.
Since $\mbox{div}(DW(\nabla\xi))$ is a lower order term in
(\ref{v_eq}), it follows that other properties of $W$ play actually no
role in the proof of our main Theorem \ref{th1} and \ref{th2}. 
We hence remark that the same results are valid when
$\mbox{div}(DW(\nabla\xi))$ is replaced by $\mbox{div}(DW((\nabla\xi)
A(X)^{-1}))$. Such term corresponds to the so-called non-Euclidean
elasticity, where the deformation $\xi$ of the reference configuration
strives to achieve a prescribed Riemannian metric $g=A^TA$ on
$\Omega$. This model pertains to the description of prestrained
materials and morphogenesis of growing tissues \cite{LePa, LMP}.

\subsection{The viscous stress tensor $\mathcal{Z}$.}
The viscous stress tensor is given by the mapping
$\mathcal{Z}:\mathbb{R}^{n\times n}\times \mathbb{R}^{n\times n}
\longrightarrow \mathbb{R}^{n\times n}$, depending on the deformation gradient $F$
and the velocity gradient $Q$.
It should be compatible with the following principles of continuum mechanics:
balance of angular momentum, frame invariance,
and the Claussius-Duhem inequality \cite{Dafbook}. That is, for every $F,Q\in
\mathbb{R}^{3\times 3}$ with $\det F > 0$, we require that:
\begin{equation}\label{visc_tensor}
\begin{minipage}{14cm}
\begin{itemize}
\item[(i)] $\mbox{skew}\left(F^{-1}\mathcal{Z}(F,Q)\right) = 0$,
i.e. $\mathcal{Z}=FS$ with $S$ symmetric.
\item[(ii)] $\mathcal{Z}(RF,R_tF + RQ) = R\mathcal{Z}(F,Q)$
for every path of rotations $R:\mathbb{R}_+\longrightarrow SO(n)$,
i.e. in view of (i): $S(RF,RKF+RQ)=S(F,Q)$ $\forall R\in SO(n)$
$\forall K\in \mbox{skew}$.
\item[(iii)] $\mathcal{Z}(F,Q):Q\geq 0$,
i.e. in view of (i): $S: \mbox{sym}(F^T Q) \geq 0$.
\end{itemize}
\end{minipage}
\end{equation}
Examples of $\mathcal{Z}$ satisfying the above are:
\begin{equation}\label{Zs}
\begin{split}
\mathcal{Z}_m(F,Q) &= [\mbox{sym}(QF^{-1})]^{2m+1} F^{-1,T},\\
\mathcal{Z}_0'(F,Q) &= 2(\mbox{det} F) \mbox{sym}(QF^{-1}) F^{-1, T},\\
\mathcal{Z}_0''(F,Q) &= 2F\mbox{sym}(F^TQ).
\end{split}
\end{equation}
We note that in the case of $\mathcal{Z}_0'$, the related Cauchy stress tensor
$T_0' = 2(\mbox{det} F)^{-1} \mathcal{Z}_2 F^T = 2\mbox{sym}(QF^{-1})$
is the Lagrangean version of the stress tensor  $2\mbox{sym}\nabla v$ written in
the Eulerian coordinates. 
For incompressible fluids $2\mbox{div}(\mbox{sym}\nabla v) = \Delta v$,
giving the usual parabolic viscous regularization of the fluid
dynamics evolutionary system.

\subsection{The main results.}
Our main assumption implying the dissipative properties of
(\ref{v_eq}) will be expressed in terms of the following condition
on a (constant coefficient) linear operator
$M:\mathbb{R}^{n\times n}\to \mathbb{R}^{n\times n}$:
\begin{equation}\label{Korni}
\|\nabla \zeta\|^2_{L_2(\mathbb{R}^n)} \leq \gamma \int_{\mathbb{R}^n} (M\nabla
\zeta):\nabla\zeta \qquad \forall \zeta\in W^1_2(\mathbb{R}^n,\mathbb{R}^n).
\end{equation}

Note that (\ref{Korni}) is a Korn-type estimate, reducing to the
classical Korn inequality for $M(F) =
\mbox{sym} F$ and $\gamma=2$ \cite{korn}.
Naturally, (\ref{Korni}) is equivalent to (\ref{Korni2}) 
which is the same estimate valid for
all $\zeta\in W_1^2(U, \mathbb{R}^n)$ with $\zeta_{|\partial U}=0$, on
any fixed open, bounded $U\subset\mathbb{R}^n$.
It can be shown, via Fourier transform (see Lemma \ref{elliptic}),
that (\ref{Korni}) is also equivalent to the strict positive definiteness
of $M$ when restricted to the space of rank-one matrices
$Q=a\otimes b$:
\begin{equation}\label{rank_one}
\forall a,b\in\mathbb{R}^n \qquad |a|^2|b|^2 = |a\otimes b|^2 \leq
\gamma \langle M(a\otimes b) : a\otimes b\rangle.
\end{equation}
We point out that the above condition resembles, naturally, the local
well-posedness criterion for the inviscid elasticity system \cite{Kato},
where the validity of (\ref{rank_one}) for $M=DW(\nabla\xi_0)$ is
equivalent to the hyperbolicity of the first-order system (\ref{v_eq})
with $\mathcal{Z}=0$.

\medskip

The main result of this paper is the following:

\begin{theorem}\label{th1}
 Let $\Omega$ be a smooth bounded domain in $\mathbb{R}^n$ 
and let:
\begin{equation}\label{reg}
\xi_0 \in W^{2}_p(\Omega), ~~~ \xi_1\in W^{2-2/p}_p(\Omega)
\quad \mbox{ for some } ~~~ p>n+2,
\end{equation}
satisfy:
$$\inf_{X\in\Omega}\det\nabla \xi_0(X) >0, \qquad \xi_0(X) = X,  ~~~ 
\xi_1(X) = 0 \quad \forall X\in\partial\Omega. $$
Assume that the viscous tensor $\mathcal{Z}$ has the property that:
\begin{equation}\label{Zgood}
\mbox{  $\forall X\in\Omega$, (\ref{Korni}) holds with  }
M=D_Q\mathcal{Z}(\nabla\xi_0(X), \nabla\xi_1(X)) 
\mbox{ and $\gamma$ independent of $X$} .
\end{equation}
Then there exists $T_{max}>0$ such that the problem (\ref{v_eq}),
(\ref{initial}), (\ref{bc}), (\ref{nonin}) admits a
unique regular solution $\xi\in W_p^{2,2}(\Omega \times (0,T))\cap
L_\infty((0,T), W_p^2(\Omega))$ 
with  $\xi_t \in W^{2,1}_p(\Omega \times (0,T))$ for all $T<T_{max}$.
\end{theorem}

The proof of Theorem \ref{th1} will be given in sections \ref{constc},
\ref{secest}, \ref{mainproof}. 
In section \ref{seclast} we show that viscous stress tensors in
(\ref{Zs}) satisfy (\ref{Zgood}): for any initial data $\xi_0, \xi_1$
in case of the linear in $Q$ tensors 
$\mathcal{Z}_0, \mathcal{Z}'_0, \mathcal{Z}''_0$, and for
initial data enjoying additionally $\det\mbox{sym}
(\nabla\xi_1(\nabla\xi_0)^{-1})\neq 0$ in case of the nonlinear (in
$Q$) tensors $\mathcal{Z}_1,\mathcal{Z}_2$, see Lemma \ref{Zexamples}.  
Thanks to this observation, Theorem \ref{th1} proves the mathematical
well-posedness of a class of physically well-posed models.

\medskip

With the same techniques of proof of
Theorem \ref{th1}, one can show that:

\begin{theorem}\label{th2}
Let $\mathcal{S}$ be the solution operator  of the problem (\ref{v_eq}) -
(\ref{nonin}) as described in Theorem \ref{th1}, given by:
$$ \mathcal{S}(\xi_0,\xi_1) = (\xi,\xi_t), $$
$$ \mathcal{S}: W^2_p(\Omega)\times W^{2-2/p}_p(\Omega)
\to \Big(W_p^{2,2}(\Omega \times (0,T))\cap L_\infty((0,T),
  W_p^2(\Omega))\Big) 
\times W^{2,1}_p(\Omega \times (0,T)).$$
Then $\mathcal{S}$ is continuous.
\end{theorem}
We omit the proof and refer instead to standard texts
\cite{Amann,LSU,Lieber,Alessandra},  or to an application of the same
methods in the more current context as in Theorem 1.2 \cite{DM} .

\subsection{Relation to previous works.}
The dynamical viscoelasticity (\ref{v_eq}) has been the subject of
vast studies in the last decades. For $\mathcal{Z}(F, Q) = Q$ 
conflicting with the frame invariance (\ref{visc_tensor}) (ii),
various results on existence, asymptotics and stability have been
obtained in \cite{1, 25, 27, 18}. For dimension $n=1$, existence of solutions
to (\ref{v_eq}) has been shown in \cite{11,4} for $\mathcal{Z}$
depending nonlinearly on $Q$.

Existence and stability of viscoelastic shock profiles for a large
class of models originating from (\ref{v_eq}) has been studied, among
others,  in \cite{AM, BLZ}.

Existence of Young measure solutions to system (\ref{v_eq}) was shown
in \cite{Demoulini}, without any additional assumptions on
$\mathcal{Z}$, but with condition (\ref{visc_tensor}) (iii)
strengthened to the uniform dissipativity i.e:
$\mathcal{Z}(F, Q) \geq \gamma |Q|^2$. 
These measure-valued solutions were shown to be the unique
classical weak solutions under the extra monotonicity assumption:
\begin{equation}\label{mon}
\langle\mathcal{Z}(F_1, Q_1) - \mathcal{Z}(F_2, Q_2) : Q_1 -
Q_2\rangle \geq \kappa |Q_1 - Q_2|^2 - l|F_1 - F_2|^2,
\end{equation}
see also \cite{Tvedt} for a treatment of slightly more general type of PDEs
under the same condition.
As noted in \cite{Demoulini}, (\ref{mon}) is incompatible with the
balance of angular momentum (\ref{visc_tensor}) (i). In particular, (\ref{mon})
is not satisfied by any of the examples in (\ref{Zs}), even
$\mathcal{Z}_0$, $\mathcal{Z}'_0$, $\mathcal{Z}''_0$ which enjoy
condition (\ref{Zgood}) for any invertible $F=\nabla\xi_0(X)$ and any
$Q=\nabla\xi_1(X)$. 

From the theory of PDEs viewpoint, our present result is a rather straightforward application of 
the theory of nonlinear (quasilinear) parabolic systems. Namely, we apply the 
maximal regularity estimates to control the nonlinearities of the system (\ref{v_eq}).
We choose the $L_p$-framework in order to avoid technical
difficulties, but a similar results and estimates are expected in
the Besov spaces framework \cite{DM}.  
In a sense, our result is hence a consequence of the classical works
of  Ladyzhenskaya, Solonnikov and Uralceva \cite{LSU}, which has been
further developed in \cite{Amann, DHP, Lieber}, and which is 
a powerful tool in the study of the parabolic-elliptic systems.

\subsection{Notation.}
By $L_p(\Omega)$ we denote the space of functions integrable with respect
to the Lebesgue measure,  with $p$-th power. 
By $W^{k,l}_p(\Omega \times (0,T))$ for $k,l \in {\mathbb N}$ we denote the anisotropic Sobolev space
defined by the norm :
\begin{equation*}
 \|u\|_{W^{k,l}_p(\Omega \times (0,T))}=\|u\|_{L_p(\Omega \times
   (0,T))} + \| \nabla^k u,\partial_t^l u\|_{L_p(\Omega \times (0,T))}, 
\end{equation*}
 where $\nabla^k$ is the $k$-th space derivative and $\partial_t$ is
 the time derivative.  The isotropic version is given by:
\begin{equation*}
 \|u\|_{W^k_p(\Omega)}=\|u\|_{L_p(\Omega)}+\|\nabla^k u\|_{L_p(\Omega)}.
\end{equation*}
The space $W^{2-2/p}_p(\Omega)$ is the trace (in time) space of
$W^{2,1}_p(\Omega \times (0,T))$. For further details we refer to \cite{BIN}.

\bigskip\bigskip

\noindent{\bf Acknowledgments.} The authors wish to thank Matthias
Hieber, Jan Pr\"uss and Vladimir Sverak for helpful consultations.
M.L. was partially supported by the NSF grant DMS-0846996,
and both authors were partially supported by the Polish MN grant N N201 547438. 

\section{The constant coefficient problem}\label{constc}
The following auxiliary result will be needed in the proof of Theorem \ref{th1}:
\begin{lemma}\label{lem1}
Assume that $M:\mathbb{R}^{n\times n}\to \mathbb{R}^{n\times n}$ is a
linear map satisfying the Korn-type inequality:
\begin{equation}\label{Korni2}
\|\nabla \zeta\|^2_{L_2(U)} \leq \gamma \int_{U} (M\nabla
\zeta):\nabla\zeta \qquad \forall \zeta\in W^1_2(U,\mathbb{R}^n)
\mbox{ with } \zeta_{|\partial U}=0
\end{equation}
Then the solution to:
\begin{equation}
\left\{\begin{array}{ll}\label{linear}
\zeta_t -\mathrm{div} \left(M \nabla \zeta\right) =f & \mbox{ in }  U \times (0,T), \\
\zeta=0 & \mbox{ on } \partial U \times (0,T),\\
\zeta(0,\cdot)=\zeta_0 & \mbox{ in } U
\end{array}\right.
\end{equation}
admits the following maximal regularity estimate:
\begin{equation}\label{para}
 \|\zeta_t, \nabla^2\zeta\|_{L_p(U \times (0,T))} \leq
 C_{\gamma,p, U}(\|f\|_{L_p(U \times (0,T))}
 +\|\zeta_0\|_{W^{2-2/p}_p(U)}),
\end{equation}
where the dependence of $C$ on $U$ is uniform for any family of domains
which are uniformly bilipschitz
homeomorphic to each other after appropriate dilations.
\end{lemma}

Towards a proof of Lemma \ref{lem1}, note first that for  $M=\mbox{Id}$,
i.e. when (\ref{Korni2}) holds trivially, (\ref{para}) is a classical
maximal regularity 
parabolic estimate for the heat equation. When $M(F) = \mbox{sym} F$,
i.e. when (\ref{Korni2}) reduces to Korn's inequality, the proof of
(\ref{para}) is also immediate. For, take $div$ of the equation in (\ref{linear}), 
and note that $\mbox{div}^T\mbox{div} (\mbox{sym}\nabla\zeta) = \mbox{div}\left(\frac{1}{2}
\Delta\zeta + \frac{1}{2}\nabla\mbox{div}\zeta\right) =
\frac{1}{2}(\mbox{div}\Delta\zeta + \frac{1}{2}\Delta\mbox{div}\zeta)
= \Delta\mbox{div}\zeta$
so that:
$$(\mbox{div}\zeta)_t - \Delta(\mbox{div}\zeta) = \mbox{div} f \quad
\mbox{ in } U\times (0,T).$$
By the maximal regularity estimate for the heat equation:
\begin{equation}\label{pomoc}
\|\nabla\mbox{div}\zeta\|_{L_p(U\times (0,T))} \leq C_{p,U} \left(
  \|f\|_{L_p(U\times (0,T))} + \|\mbox{div} \zeta_0\|_{W_p^{1-1/p}(U)}\right).
\end{equation}
Now, (\ref{linear}) can be written as:
$$\zeta_t - \frac{1}{2} \Delta\zeta = f + \frac{1}{2}
\nabla\mbox{div}\zeta \quad \mbox{ in } U\times (0,T).$$
We hence obtain:
\begin{equation*}
\|\zeta_t, \nabla^2\zeta\|_{L_p(U\times (0,T))} \leq C_{p,U} \left(
  \|f + \frac{1}{2}\nabla\mbox{div}\zeta \|_{L_p(U\times (0,T))} 
+ \|\mbox{div} \zeta_0\|_{W_p^{2-1/p}(U)}\right),
\end{equation*}
which combined with (\ref{pomoc}) yields (\ref{para}).

\bigskip

In the general case, Lemma \ref{lem1} follows from the maximal
regularity theory developed for parabolic initial-boundary value problems 
in \cite{DHP}. Under the ellipticity condition (b) on page 98 in there (see also
Definition 5.1),  the estimate (\ref{para}) is a consequence of 
Theorem 7.11.  We now prove that condition (\ref{Korni2}) implies that
the constant coefficient operator $- {\rm div\,} (M\nabla \zeta)$ has
its spectrum contained in the proper sector of the complex plane, 
which immediately gives ellipticity in the sense of \cite{DHP}. 

\begin{lemma}\label{elliptic}
Conditions (\ref{Korni}), (\ref{Korni2}) and (\ref{rank_one}) are equivalent.
Moreover, under any of these conditions the operator  $-{\rm div\,}
M\nabla (\cdot) $  is elliptic, i.e:
\begin{equation}\label{x3}
 {\rm spec}\big(-{\rm div\,} M\nabla (\cdot)\big) \subset \Big\{ z\in
 \mathbb{C} : \mathrm{Re}~ z > 0, \mbox{ and } \mathrm{arg}~  z 
 \in [\alpha_*,\alpha^*] \mbox{ with  }
 -\frac{\pi}{2}< \alpha_* < \alpha^* <\frac{\pi}{2}\Big\}.
\end{equation} 
\end{lemma}
\begin{proof}
{\bf 1.} Conditions (\ref{Korni}) and (\ref{Korni2}) are equivalent in
view of the density of $C_c^\infty(\mathbb{R}^n)$ in
$W_1^2(\mathbb{R}^n)$. To include (\ref{rank_one}), we use linearity
of Fourier transform and Plancherel's identity:
\begin{equation*}
\begin{split}
& \|\nabla\zeta\|_{L_2(\mathbb{R}^n})^2 =
\|(\nabla\zeta)^\wedge \|_{L_2(\mathbb{R}^n})^2 = \int_{\mathbb{R}^n}
|\hat\zeta(k)\otimes k|^2~\mbox{d}k\\
& \int_{\mathbb{R}^n} \langle M(\nabla\zeta):\nabla\zeta\rangle
= \int_{\mathbb{R}^n} \langle M((\nabla\zeta)^\wedge): \overline{(\nabla\zeta)^\wedge}\rangle
= \int_{\mathbb{R}^n} \langle M(\hat\zeta(k)\otimes k) :
\overline{(\hat\zeta(k)\otimes k)}\rangle ~ \mbox{d}k.
\end{split}
\end{equation*}
Hence, (\ref{Korni}) is equivalent to:
\begin{equation}\label{x1}
\forall \zeta\in W_1^2(\mathbb{R}^n) \qquad \|\hat\zeta \otimes
k\|_{L_2(\mathbb{R}^n)}^2 \leq \gamma \int_{\mathbb{R}^n} \left\langle M(\hat\zeta\otimes k) :
\overline{(\hat\zeta\otimes k)}~ \right\rangle ~ \mbox{d}k.
\end{equation}
It is therefore clear that (\ref{rank_one}) implies (\ref{Korni}). On
the other hand, given $k_0, a_0\in\mathbb{R}^n$, consider:
$\hat\zeta_m(k) = \left(\rho_m^{1/2}(k-k_0) + \rho_m^{1/2}(k+k_0)\right) a_0$,
where $\rho_m$ is the standard radially symmetric mollifier supported
in the ball $B(0, 1/m)$. Applying (\ref{x1}) to $\xi_n\in
W_1^2(\mathbb{R}^n, \mathbb{R}^n)$ and passing to the limit
$m\to\infty$, yield (\ref{rank_one}) for the matrix $Q=a_0\otimes k_0$.

\medskip

{\bf 2.} To prove (\ref{x3}), consider the eigenvalue problem:
$$\lambda\zeta - \mbox{div}(M(\nabla\zeta)) = 0 \qquad \mbox{in } \mathbb{R}^n,$$
which after passing to the Fourier variable $k\in \mathbb{R}^n$ becomes:
\begin{equation}\label{x8}
\lambda\hat\zeta(k) = M(\hat\zeta (k) \otimes k) k. 
\end{equation}
Upon writing $\lambda = \sigma |k|^2$, the problem (\ref{x8}) is
equivalent to locating the eigenvalues $\sigma$ of the family of linear
operators $\{M_k\}_{|k|=1}$, $M_k:\mathbb{R}^n\to \mathbb{R}^n$ given by:
$M_k(a) = M(a\otimes k)k$. Recalling (\ref{rank_one}) we see that each
$M_k$ is strictly positive definite:
$$ M_k(a) \cdot a = \big\langle M(a\otimes k) : (a\otimes k)\big\rangle 
\geq\frac{ |k|^2}{\gamma} |a|^2.$$
Consequently,  spectrum of every $M_k$ satisfies $\mbox{Re}~\sigma > 0$.
By continuity with respect to $k$ which varies in the compact set
$|k|=1$, we obtain the inclusion (\ref{x3}).
\end{proof}

\bigskip

Finally, we have the following:

\begin{lemma}\label{Zexamples}
The viscous stress tensors $\mathcal{Z}$ in (\ref{Zs}) satisfy
(\ref{Korni2}) with $M=D_Q\mathcal{Z}(F_0, Q_0)$,  
for every $F_0, Q_0$ with $\det F_0 > 0$,
in the following manner:
\begin{itemize}
\item[(i)] $\mathcal{Z}''_0$  with  $\gamma = |F_0^{-1,T}|^2$.
\item[(ii)] $\mathcal{Z}'_0$ with  $\gamma = |F_0|^2 (\det F_0)^{-1}$.
\item[(iii)] $\mathcal{Z}_0$ with $\gamma = \frac{1}{2}|F_0|^2$.
\end{itemize}
If we additionally assume that $\det\mathrm{sym}(Q_0F_0^{-1})\neq 0$
then we also have:
\begin{itemize} 
\item[(iv)] $\mathcal{Z}_1$ with $\gamma = 2|F_0|^2|\mathrm{sym}
  (Q_0F_0^{-1})^{-1}|^2$.
\item[(v)] $\mathcal{Z}_2$ with  $\gamma = 2|F_0|^2|\mathrm{sym}
  (Q_0F_0^{-1})^{-1}|^4$.
\end{itemize}
\end{lemma}
The proof of Lemma \ref{Zexamples} will be given in section \ref{seclast}.
We now remark that in the proof of the main Theorem \ref{th1}, Lemma \ref{lem1}
will be used to the operators
$M=M_X=D_Q\mathcal{Z}(F_{0}, Q_0)$, at finitely many  
spacial points $X\in\Omega$,
where $F_0=\nabla\xi_0(X)$ and 
$Q_0=\nabla\xi_1(X)$. It is clear that when the initial
data $\xi_0$, $\xi_1$ with regularity (\ref{reg}) satisfy  $\det\nabla
\xi > 0$ (or the two conditions $\det\nabla\xi > 0$ and
$\det\mbox{sym} (\nabla\xi_1 (\nabla\xi_0)^{-1})\neq 0$ whenever
required) then the constants
$\gamma$ in Lemma \ref{Zexamples} have a common upper and lower bounds,
independent of $X$. Therefore, Lemma \ref{lem1} and the estimate
(\ref{para}) may be used with a uniform constant $C_{p,U}$, also
independent of $X$.

\section{The main a-priori estimate}\label{secest}

Given $\xi_0$, $\xi_1$ as in Theorem \ref{th1}, let 
$\bar\xi_1\in W_p^{2,1}(\Omega\times \mathbb{R}_+)$ be
the solution to:
\begin{equation*}
\left\{\begin{array}{ll} 
(\bar\xi_1)_t - \Delta \bar\xi_1 = 0  & \mbox{ in }  \Omega \times \mathbb{R}_+, \\
\bar\xi_1=0 & \mbox{ on } \partial\Omega \times \mathbb{R}_+,\\
\bar\xi_1(0,\cdot)=\xi_1 & \mbox{ in } \Omega,
\end{array}\right.
\end{equation*}
Define the extension $\bar\xi$ of $\xi_0$, so that $\partial_t\bar\xi =\bar\xi_1$:
\begin{equation}\label{barxi}
\bar\xi(t,x) = \xi_0(x) + \int_0^t \bar\xi_1(s,x)~\mbox{d}s.
\end{equation}
By continuity, it is clear that: $\inf_{\Omega\times(0,T)}\det\nabla\xi > 0$ for $T$ sufficiently small.
We define:
$$D=D(T)= \|\bar\xi_{tt}, \nabla^2\bar\xi_t\|_{L_p(\Omega\times (0,T))}$$
and note that:
\begin{equation}\label{Dsmall}
\lim_{T\to 0} D(T) = 0.
\end{equation}

\begin{lemma}\label{lemik}
Let $\xi_0$, $\xi_1$ be as in Theorem \ref{th1} and assume that:
\begin{equation*}
\mbox{ for every $X\in\Omega$ (\ref{Korni}) holds with  }
M=D_Q\mathcal{Z}(\nabla\xi_0(X), \nabla\xi_1(X)).
\end{equation*}
Let $\xi\in W_p^{2,2}(\Omega \times (0,T_0))$ 
with  $\xi_t \in W^{2,1}_p(\Omega \times (0,T_0))$ be a solution to the problem (\ref{v_eq}),
(\ref{initial}), (\ref{bc}), (\ref{nonin}), and denote:
$$\Theta=\Theta(T)= \|(\xi-\bar\xi)_{tt},
\nabla^2(\xi-\bar\xi)_t\|_{L_p(\Omega\times (0,T))},$$
where $\bar\xi$ is as in (\ref{barxi}).
Then, there exists $T_{00}<T_0$ and a constant $C$, both depending
only on $\xi_0$ and $\xi_1$ (and, naturally, on $\Omega$ and $p$), such that
for every $T<T_{00}$ we have:
\begin{equation}\label{ineqlem}
\Theta\leq C (T^{1/p} + D + (T^{1/p} + D)\Theta + \Theta^2 + \Theta^4).
\end{equation}
In particular:
\begin{equation}\label{bdlem}
\Theta(T)\leq C \qquad \forall T< T_{00}.
\end{equation}
\end{lemma}

Before we give the proof of the lemma, we gather below some standard inequalities
that will be frequently used for different functions:
$u$ defined on $\Omega\times (0,T)$, and $w$ defined on $\Omega$.
We always assume that $T<1$.
\begin{equation}\label{ell}
\|w\|_{W^2_p(\Omega)} \leq C_{p,\Omega} \|\Delta w\|_{L_p(\Omega)}
\quad \mbox{ when } w_{|\partial\Omega} = 0,
\end{equation}
\begin{equation}\label{par}
\begin{split}
\sup_{t\in(0,T)}\|u(t,\cdot)\|_{W^{2-2/p}_p(\Omega)} \leq
C_{p,\Omega} &\left(\|u_t - \Delta u\|_{L_p(\Omega\times (0,T))} + \|u(0,\cdot)\|_{W_p^{2-2/p}(\Omega)}\right)\\
 &\qquad\qquad\qquad\qquad\qquad \mbox{ when } u_{|\partial\Omega\times (0,T)} = 0,
\end{split}
\end{equation}
\begin{equation}\label{em1}
\|w\|_{L_\infty(\Omega)} \leq C_{p,\Omega}
\|w\|_{W_p^{1-1/p}(\Omega)}, ~~~\quad \mbox{ in fact: } ~~~
\|w\|_{\mathcal{C}^{0,\alpha}(\Omega)}\leq C_{\alpha,p,\Omega}\|w\|_{W_p^{1-1/p}(\Omega)},
\end{equation}
\begin{equation}\label{em2}
\|\nabla u\|_{\mathcal{C}^{\alpha,\alpha/2}}\leq C_{\alpha,p,U}\|u\|_{W_p^{2,1}(U\times (0,T))}.
\end{equation}
The inequality (\ref{ell}) is the usual elliptic estimate \cite{GT},
and (\ref{par}) is the parabolic estimate from \cite{BIN}. The
Morrey embedding gives (\ref{em1}) for $p>n+2$ \cite{GT}, while
(\ref{em2}) follows from the embedding $\nabla W_p^{2,1} (\Omega\times
(0,T))\subset L_\infty(\Omega \times (0,T))$, also valid for $p>n+2$ \cite{LSU}.
We stress that the constants $C$ in all the above bounds are
universal, i.e. they are independent of $T$.
Additionally, the dependence of $C$ in (\ref{em2}) on $U\subset\Omega$
is uniform for any family of domains
which are uniformly bilipschitz homeomorphic to each other after appropriate dilations.


\smallskip

We further remark the following simple bound:
\begin{equation}\label{a}
\begin{split}
\|\nabla^2u&\|_{L_p(\Omega\times(0,T))}  = \left(\int_\Omega \int_0^T \left|\int_0^t
\nabla^2u_t ~\mbox{d}s + \nabla^2u(0,\cdot)\right|^p ~\mbox{d}t ~\mbox{d}X\right)^{1/p}\\
& \leq  T^{1/p} \left(\int_\Omega T^{p/p'} \int_0^T|\nabla^2
  u_t|^p ~\mbox{d}t ~\mbox{d}X\right)^{1/p} +
T^{1/p}\|\nabla^2u(0,\cdot)\|_{L_p(\Omega)}\\ 
& = T^{1/p}\left(\|\nabla^2u_t\|_{L_p(\Omega\times (0,T))} + \|\nabla^2u(0,\cdot)\|_{L_p(\Omega)}\right).
\end{split}
\end{equation}

\smallskip

Let now $\xi$ and $\bar\xi$ be as in Lemma \ref{lemik}. 
Using (\ref{par}) to $(\xi-\bar\xi)_t$ we obtain:
\begin{equation}\label{b}
\begin{split}
\sup_{t\in (0,T)} &\|(\xi-\bar\xi)_t(t,\cdot)\|_{L_p(\Omega)} + 
\sup_{t\in (0,T)} \|\nabla(\xi-\bar\xi)_t(t,\cdot)\|_{L_p(\Omega)}  \\
& \leq
C \left(\|(\xi-\bar\xi)_{tt}\|_{L_p(\Omega\times (0,T))} +
  \|\nabla^2(\xi-\bar\xi)_t \|_{L_p(\Omega\times (0,T))} \right) \leq C\Theta,
\end{split}
\end{equation}
and consequently:
\begin{equation}\label{b1}
\|(\xi-\bar\xi)_t \|_{L_p(\Omega\times (0,T))} +
\|\nabla(\xi-\bar\xi)_t \|_{L_p(\Omega\times (0,T))} 
\leq C T^{1/p} \Theta.
\end{equation}
By (\ref{em1}), (\ref{par}) used to $\xi-\bar\xi$, and (\ref{a}),
(\ref{b1}) we get:
\begin{equation}\label{b2}
\begin{split}
\|\nabla(\xi-\bar\xi) &\|_{L_\infty(\Omega\times (0,T))}  
= \sup_{t\in (0,T)}\|\nabla(\xi-\bar\xi)(t,\cdot) \|_{L_\infty(\Omega)} \\ & \leq
\sup_{t\in (0,T)} \|\nabla(\xi-\bar\xi)(t,\cdot)
\|_{W_p^{1-1/p}(\Omega)} \\ & \leq C\left( \|(\xi-\bar\xi)_{t}\|_{L_p(\Omega\times (0,T))} +
  \|\nabla^2(\xi-\bar\xi) \|_{L_p(\Omega\times (0,T))} \right) \leq C T^{1/p}\Theta. 
\end{split}
\end{equation}
Likewise, using (\ref{em1}) and (\ref{par}) to $(\xi-\bar\xi)_t$, we
directly obtain:
\begin{equation}\label{b3}
\|(\xi-\bar\xi)_t\|_{L_\infty(\Omega\times (0,T))}  + 
\|\nabla(\xi-\bar\xi)_t\|_{L_\infty(\Omega\times (0,T))}  \leq C \Theta. 
\end{equation}

In all the above inequalities (\ref{b}) -- (\ref{b3}), we write
$\Theta=\Theta(T)$. The constant $C$ depends only on the initial data
of the problem $\xi_0, \xi_1$ (in addition to its dependence on $\Omega$ and $p$).

\bigskip\bigskip

\noindent {\bf Proof of Lemma \ref{lemik}.}

\noindent 
We will always assume that $T<1$. Note that for $T<T_0$ sufficiently small, the constraint
(\ref{nonin}) is a consequence of the same constraint on the initial
data $\xi_0$, by continuity. Likewise:
\begin{equation}\label{global}
\|D\mathcal{Z}(\nabla\bar\xi,\nabla\bar\xi_t), D^2\mathcal{Z}(\nabla\bar\xi,\nabla\bar\xi_t),
D^3\mathcal{Z}(\nabla\bar\xi,\nabla\bar\xi_t)\|_{L_\infty(\Omega\times(0,T))} \leq C.
\end{equation}

\smallskip

{\bf 1.}  The system (\ref{v_eq}) can be rewritten as:
\begin{equation*}
(\xi-\bar\xi)_{tt} - \mbox{div}
\left(\mathcal{Z}(\nabla\xi,\nabla\xi_t) - \mathcal{Z}(\nabla\bar\xi,
    \nabla\bar\xi_t)\right) = \mbox{div} \left(DW(\nabla \xi)\right) + 
\mbox{div}\left(\mathcal{Z}(\nabla\bar\xi,\nabla\bar\xi_t)\right) - \bar\xi_{tt} 
\end{equation*}
and further, it has the form:
\begin{equation}\label{v_eq2}
(\xi-\bar\xi)_{tt} - \mbox{div}
\left(D_Q\mathcal{Z}(\nabla\bar\xi,\nabla\bar\xi_t)\nabla(\xi-\bar\xi)_t\right)
= F[\xi,\bar\xi],
\end{equation}
where:
\begin{equation}\label{Feq}
\begin{split}
&F[\xi,\bar\xi]  =  \mbox{div} \left(DW(\nabla \xi)\right) + 
\mbox{div}\left(\mathcal{Z}(\nabla\bar\xi,\nabla\bar\xi_t)\right) -
\bar\xi_{tt} \\ & \qquad \qquad + \mbox{div}
\left(\mathcal{Z}(\nabla\xi,\nabla\xi_t) - \mathcal{Z}(\nabla\bar\xi,
    \nabla\xi_t)\right) \\ & \qquad \qquad + \mbox{div}
\left(\mathcal{Z}(\nabla\bar\xi,\nabla\xi_t) - \mathcal{Z}(\nabla\bar\xi,
    \nabla\bar\xi_t) -
    D_Q\mathcal{Z}(\nabla\bar\xi,\nabla\bar\xi_t)\nabla(\xi-\bar\xi)_t\right) \\
& = \mbox{div} \left(DW(\nabla \xi)\right) + 
\mbox{div}\left(\mathcal{Z}(\nabla\bar\xi,\nabla\bar\xi_t)\right) -
\bar\xi_{tt} \\ & \quad + \mbox{div}\left(D_F\mathcal{Z}(\nabla\bar\xi,\nabla\xi_t)
  \nabla(\xi-\bar\xi)\right)  \\ & \quad 
+ \mbox{div}\left(\int_0^1 (1-s) D^2_{FF}\mathcal{Z}(s\nabla\xi +
  (1-s) \nabla\bar\xi,\nabla\xi_t)  (\nabla(\xi-\bar\xi)\otimes
  \nabla(\xi-\bar\xi))~\mbox{d}s\right) \\ & \quad 
+ \mbox{div}\left(\int_0^1 (1-s) D^2_{QQ}\mathcal{Z}(\nabla\bar\xi, s\nabla\xi_t +
  (1-s) \nabla\bar\xi_t)  (\nabla(\xi-\bar\xi)_t\otimes \nabla(\xi-\bar\xi)_t)~\mbox{d}s\right).
\end{split}
\end{equation}
We shall now prove the bound:
\begin{equation}\label{firstbound}
\|F[\xi,\bar\xi]\|_{L_p(\Omega\times (0,T)} \leq C\left(T^{1/p} + D +
  (T^{1/p} + D)\Theta + \Theta^2 + \Theta^4\right).
\end{equation}

By (\ref{a}) and (\ref{b2}) it follows that:
\begin{equation}\label{s1}
\begin{split}
\| \mbox{div} & \left(DW(\nabla \xi)\right)\|_{L_p(\Omega\times (0,T))} 
\leq \| D^2W(\nabla \xi)\|_{L_\infty(\Omega\times (0,T))}
\|\nabla^2\xi\|_{L_p(\Omega\times (0,T))}  \\ & \leq \left(\| D^2W(\nabla \bar\xi)\|_{L_\infty}
+ C\|\nabla(\xi-\bar\xi)\|_{L_\infty}
\right) T^{1/p} \cdot \\
& \qquad\qquad\qquad \qquad\qquad 
\cdot \left(\|\nabla^2(\xi-\bar\xi)_t\|_{L_p} + 
  \|\nabla^2\bar\xi_t\|_{L_p}  + \|\nabla^2\xi_0\|_{L_p(\Omega)}
\right) \\ & \leq 
C(1+T^{1/p}\Theta) T^{1/p}(1+ \Theta + D)\leq C
T^{1/p}(1+\Theta)(1+\Theta + D).
\end{split}
\end{equation}
Using (\ref{global}) and (\ref{a}) to $\bar\xi$, we obtain:
\begin{equation}\label{s2}
\begin{split}
\| \mbox{div} & \left(\mathcal{Z}(\nabla \bar\xi,
  \nabla\bar\xi_t)\right)\|_{L_p(\Omega\times (0,T))} \\ & 
\leq \| D\mathcal{Z}(\nabla \bar\xi, \nabla\bar\xi_t)\|_{L_\infty(\Omega\times (0,T))}
\left(\|\nabla^2\bar\xi\|_{L_p(\Omega\times (0,T))} +
  \|\nabla^2\bar\xi_t\|_{L_p(\Omega\times (0,T))} \right) \\ & 
\leq C \left( T^{1/p} \|\nabla^2\bar\xi_t\|_{L_p} +  T^{1/p}
  \|\nabla^2\xi_0\|_{L_p(\Omega)} +  \|\nabla^2\bar\xi_t\|_{L_p}\right) 
\leq C(T^{1/p}+ D).
\end{split}
\end{equation}
By (\ref{global}), (\ref{b3}), (\ref{a}), (\ref{b2}) we get:
\begin{equation}\label{s3}
\begin{split}
& \| \mbox{div}  \left(D_F\mathcal{Z}(\nabla \bar\xi,
  \nabla\xi_t)\nabla(\xi-\bar\xi)\right)\|_{L_p(\Omega\times (0,T))} \\ & 
\leq \| D^2_{FF}\mathcal{Z}(\nabla \bar\xi, \nabla\xi_t)\|_{L_\infty}
\left(\|\nabla^2\bar\xi\|_{L_p(\Omega\times (0,T))} +
  \|\nabla^2\bar\xi_t\|_{L_p(\Omega\times (0,T))} \right) 
\|\nabla(\xi-\bar\xi)\|_{L_\infty} \\ &
\qquad\qquad\qquad \qquad\qquad \qquad\qquad 
 + \| D_{F}\mathcal{Z}(\nabla \bar\xi, \nabla\xi_t)\|_{L_\infty}
\|\nabla^2(\xi-\bar\xi)\|_{L_p(\Omega\times (0,T))}
\\ &  \leq C\left(1+ \|\nabla(\xi-\bar\xi)_t\|_{L_\infty}\right)
\left(\|\nabla^2\bar\xi\|_{L_p} + \|\nabla^2\bar\xi_t\|_{L_p} + \|\nabla^2(\xi-\bar\xi)_t\|_{L_p} \right)
\|\nabla(\xi-\bar\xi)\|_{L_\infty} \\ &
\qquad\qquad\qquad \qquad\qquad \qquad\qquad 
+ C \left(1+ \|\nabla(\xi-\bar\xi)_t\|_{L_\infty}\right)
\|\nabla^2(\xi-\bar\xi)\|_{L_p} \\ &
\leq C (1+\Theta)( T^{1/p} + \Theta + D) T^{1/p}\Theta + C
  (1+\Theta) T^{1/p}\Theta \\ & \leq CT^{1/p}(1+\Theta+ D)\Theta (1+\Theta).
\end{split}
\end{equation}
and:
\begin{equation}\label{s4}
\begin{split}
& \| \mbox{div}  \left(\int_0^1(1-s) D^2_{FF}\mathcal{Z}(s\nabla\xi + (1-s)\nabla\bar\xi,\nabla\xi_t)
(\nabla(\xi-\bar\xi)\otimes\nabla(\xi-\bar\xi)) ~\mbox{d}s\right)\|_{L_p(\Omega\times (0,T))} \\ & 
\leq \sup_{s\in [0,1]} 
\| \mbox{div}  \left(D^2_{FF}\mathcal{Z}(s\nabla\xi + (1-s)\nabla\bar\xi,\nabla\xi_t)
(\nabla(\xi-\bar\xi)\otimes\nabla(\xi-\bar\xi))
\right)\|_{L_p(\Omega\times (0,T))} \\ & 
\leq \sup_{s\in [0,1]} \Bigg[\| D^3\mathcal{Z}(s\nabla\xi +(1-s)\nabla\bar\xi, \nabla\xi_t)\|_{L_\infty}
\left(\|\nabla^2\bar\xi\|_{L_p} + \|\nabla^2(\xi-\bar\xi)\|_{L_p} +
  \|\nabla^2\xi_t\|_{L_p} \right)  \|\nabla(\xi-\bar\xi)\|_{L_\infty}^2 \\ &
\qquad\qquad\qquad \qquad\qquad
 + \| D^2\mathcal{Z}(s\nabla\xi + (1-s)\nabla\bar\xi, \nabla\xi_t)\|_{L_\infty}
\|\nabla^2(\xi-\bar\xi)\|_{L_p} \|\nabla(\xi-\bar\xi)\|_{L_\infty}
\Bigg]
\\ &  \leq C\left(1+ \|\nabla(\xi-\bar\xi)\|_{L_\infty} +
  \|\nabla(\xi-\bar\xi)_t\|_{L_\infty}\right) (T^{1/p} + \Theta + D)
T^{2/p}\Theta^2 \\ & 
\qquad\qquad\qquad \qquad\qquad
+ C \left(1+ \|\nabla(\xi-\bar\xi)\|_{L_\infty} + \|\nabla(\xi-\bar\xi)_t\|_{L_\infty}\right)
 T^{2/p}\Theta^2 \\ & \leq CT^{1/p}(1+\Theta+ D)\Theta^2 (1+\Theta).
\end{split}
\end{equation}
In the same manner, we see that:
\begin{equation}\label{s5}
\begin{split}
\| \mbox{div}  &\left(\int_0^1(1-s)
  D^2_{QQ}\mathcal{Z}(\nabla\bar\xi, s\nabla\xi_t + (1-s)\nabla\bar\xi_t)
(\nabla(\xi-\bar\xi)_t\otimes\nabla(\xi-\bar\xi)_t) ~\mbox{d}s\right)\|_{L_p(\Omega\times (0,T))} \\ & 
\leq C\left(1+ \|\nabla(\xi-\bar\xi)_t\|_{L_\infty} \right)
\left(\|\nabla^2\bar\xi\|_{L_p} + \|\nabla^2(\xi-\bar\xi)_t\|_{L_p} +
  \|\nabla^2\xi_t\|_{L_p} \right)  \|\nabla(\xi-\bar\xi)_t\|_{L_\infty}^2 \\ &
\qquad\qquad\qquad \qquad\qquad
 + C \left(1+ \|\nabla(\xi-\bar\xi)_t\|_{L_\infty} \right)
\|\nabla^2(\xi-\bar\xi)_t\|_{L_p} \|\nabla(\xi-\bar\xi)_t\|_{L_\infty}
\\ &  \leq C\left(1+\Theta\right) (T^{1/p} + \Theta + D)\Theta^2 
+ C \left(1+ \Theta \right)\Theta\\ & \leq C(T^{1/p}+\Theta+ D)\Theta (1+\Theta)^2.
\end{split}
\end{equation}
Combining (\ref{s1}) -- (\ref{s5}), the bound (\ref{firstbound})
follows if only $T<1$, ensuring $D(T) < 1$ by (\ref{Dsmall}).

\bigskip

{\bf 2.} We will now work with the localizations of the system (\ref{v_eq2}). Let
$\{B_k\}_{k=1}^{N}$ be a covering of $\Omega$ by a finite number
$N=N(r)$ of balls $B_k = B(X_k,r)$ 
with centers $X_k\in\Omega$ and radius $r<1$. This family of
coverings (parametrized by $r$) should such that 
all the sets $2B_k\cap\Omega$ are uniformly bilipschitz homeomorphic  
to each other after appropriate dilations and that the covering
numbers of $\{2B_k\cap\Omega\}_k$ are independent of $r$. 

Let $\pi_k:\mathbb{R}^n\to [0,1]$ be
smooth cut-off functions satisfying: $\pi_k = 1$ on $B_k$, and $\pi_k =
0$ on $\mathbb{R}^n\setminus 2B_k$ where $2B_k = B(X_k, 2r)$, and
$\|\nabla^\alpha\pi_k\|_{L_\infty} \leq C r^{-|\alpha|}$. After
multiplying (\ref{v_eq2}) by $\pi_k$, we obtain:
\begin{equation} \label{v_eq3k}
\big(\pi_k(\xi-\bar\xi)\big)_{tt} - \mbox{div} \Big( D_Q\mathcal{Z}
  (\nabla\xi_0(X_k), \nabla\xi_1(X_k))\nabla \big(\pi_k(\xi-\bar\xi)_t\big)\Big)
= \pi_k F[\xi,\bar\xi] + G_k[\xi,\bar\xi], 
\end{equation}
where:
\begin{equation*}
\begin{split}
G_k[\xi,\bar\xi]  = & \pi_k\mbox{ div} \Big([D_Q\mathcal{Z}(\nabla\bar\xi,
\nabla\bar\xi_t) - D_Q\mathcal{Z}(\nabla\xi_0(X_k),
\nabla\xi_1(X_k))]\nabla (\xi-\bar\xi)_t\Big)\\ & 
- \Big(D_Q\mathcal{Z}(\nabla\xi_0(X_k), \nabla\xi_1(X_k)) \nabla(\xi-
  \bar\xi)_t\Big) \nabla\pi_k \\ & -
  \mbox{div}\left(D_Q\mathcal{Z}(\nabla\xi_0(X_k), \nabla\xi_1(X_k) )
      ((\xi-\bar\xi)_t\otimes \nabla\pi_k)\right).
\end{split}
\end{equation*}
We shall now prove the bound:
\begin{equation}\label{secondbound}
\begin{split}
\|G_k[\xi, \bar\xi]\|_{L_p(2B_k\times (0,T))} \leq &
C(r^\alpha + T^{\alpha/2})
\|\pi_k\nabla^2(\xi-\bar\xi)_t\|_{L_p(2B_k\times (0,T))} \\ & + 
C(1+\frac{1}{r^2})
(T^{1/p} + D) \Theta (1+\Theta).
\end{split}
\end{equation}

Using (\ref{b1}), we obtain:
\begin{equation}\label{p2}
\begin{split}
\|&\left(D_Q\mathcal{Z}(\nabla\xi_0(X_k), \nabla\xi_1(X_k))
\nabla(\xi-\bar\xi)_t\right)\nabla\pi_k\|_{L_p(2B_k\times(0,T))} \\ & 
\qquad \qquad\qquad\qquad\qquad\qquad\qquad\qquad
\leq \frac{C}{r} \|\nabla(\xi-\bar\xi)_t\|_{L_p} \leq \frac{C}{r} T^{1/p}\Theta^2. 
\end{split}
\end{equation}
Likewise:
\begin{equation}\label{p3}
\begin{split}
\|\mbox{div}&\left(D_Q\mathcal{Z}(\nabla\xi_0(X_k), \nabla\xi_1(X_k))
((\xi-\bar\xi)_t\otimes \nabla\pi_k)\right)\|_{L_p(2B_k\times(0,T))}
\\ &  \qquad\qquad\qquad\qquad\qquad
\frac{C}{r} \|\nabla(\xi-\bar\xi)_t\|_{L_p} +  \frac{C}{r^2}\|(\xi-\bar\xi)_t\|_{L_p}
\leq \frac{C}{r^2} T^{1/p}\Theta. 
\end{split}
\end{equation}
Finally, by (\ref{global}), (\ref{a}), (\ref{b3}) and (\ref{em2}) we have:
\begin{equation}\label{p1}
\begin{split}
& \|\pi_k \mbox{ div} \Big([D_Q\mathcal{Z}(\nabla\bar\xi,
\nabla\bar\xi_t) - D_Q\mathcal{Z}(\nabla\xi_0(X_k),
\nabla\xi_1(X_k))]\nabla (\xi-\bar\xi)_t\Big)\|_{L_p(2B_k\times(0,T))} 
\\ & \leq C\left(\|\nabla^2\bar\xi\|_{L_p} +  \|\nabla^2\bar\xi_t\|_{L_p}\right)
\|\nabla(\xi-\bar\xi)_t\|_{L_\infty} \\ &
\qquad + \|D_Q\mathcal{Z}(\nabla\bar\xi, \nabla\bar\xi_t) - D_Q\mathcal{Z}(\nabla\xi_0(X_k),
\nabla\xi_1(X_k))\|_{L_\infty(2B_k\times (0,T))} \|\pi_k\nabla^2
(\xi-\bar\xi)_t \|_{L_p(2B_k\times (0,T))}\\ & 
\leq  C (T^{1/p} + D)\Theta  \\ & 
\qquad + C\left( \|\nabla\bar\xi-\nabla\xi_0(X_k)\|_{L_\infty(2B_k\times(0,T))}
  + \|\nabla\bar\xi_t -
  \nabla\xi_1(X_k)\|_{L_\infty(2B_k\times(0,T))}\right)\cdot \\ & 
\qquad\qquad\qquad\qquad\qquad\qquad\qquad\qquad\qquad\qquad\qquad\qquad
\cdot \|\pi_k\nabla^2 (\xi-\bar\xi)_t \|_{L_p(2B_k\times (0,T))} \\ &
\leq C (T^{1/p} + D)\Theta 
+ C (r^\alpha + T^{\alpha/2}) \|\bar\xi\|_{W_p^{2,1}}
\|\pi_k\nabla^2(\xi-\bar\xi)_t\|_{L_p(2B_k\times (0,T))}.
\end{split}
\end{equation}
Combining (\ref{p3}) -- (\ref{p1}) and noting that 
$\|\bar\xi\|_{W_p^{2,1}(\Omega\times (0,T))}\leq
C(\|\xi_0\|_{W^2_p(\Omega)} +\|\xi_1\|_{W^{2-2/p}_p(\Omega)})$, 
for $T$ small enough we conclude (\ref{secondbound}) in view of (\ref{Dsmall}).

\medskip

We now use Lemma \ref{lem1} to the problem (\ref{v_eq3k}), i.e. we set
$\zeta= \pi_k(\xi-\bar\xi)_t$, $M=D_Q\mathcal{Z}(\nabla\xi_0(X_k),
\nabla\xi_1(X_k))$ where $U=2B_k\cap \Omega$ is the uniform constant
from the assumption (\ref{Zgood}). Indeed, it is easy to notice that
if (\ref{Korni2}) holds for some set $U$ then it holds with the same
constant $\gamma$ on every open subset $U_1\subset U$. By (\ref{para})
we now obtain:
\begin{equation*}
\|\pi_k(\xi-\bar\xi)_{tt},
\nabla^2(\pi_k(\xi-\bar\xi)_t)\|_{L_p(2B_k\times (0,T))} \leq 
C \|\pi_k F[\xi, \bar\xi], G_k[\xi, \bar\xi]\|_{L_p(2B_k\times (0,T))}. 
\end{equation*}
Summing over finitely many $k:1\ldots N$, we get in view of (\ref{secondbound}):
\begin{equation*}
\begin{split}
\|(\xi-\bar\xi)_{tt},&
\nabla^2(\xi-\bar\xi)_t\|_{L_p(\Omega\times (0,T))}  \leq 
C (r^\alpha + T^{\alpha/2})
\|\nabla^2(\xi-\bar\xi)_t\|_{L_p(\Omega\times (0,T))} \\ & +
CN^{1/p}\Big((1+1/r)(T^{1/p} + D)\Theta (1+\Theta)
+ \|F[\xi, \bar\xi]\|_{L_p(\Omega\times (0,T))}\Big), 
\end{split}
\end{equation*}
where, again, $C$ depends only on the covering number of
$\{B_k\}_{k=1}^N$, on $\bar\xi$, $p$ and $\Omega$, but not on $r, N,
T$ or $\Theta$. Consequently, for $r$ and $T$ sufficiently small, we
arrive at:
$$\Theta \leq CN^{1/p} \left(1+\frac{1}{r}\right)\big(T^{1/p} +D +
(T^{1/p} + D) \Theta +\Theta^2 + \Theta^4\big) $$
in virtue of (\ref{firstbound}).
This concludes the proof of (\ref{ineqlem}).

\bigskip

{\bf 3.} To prove (\ref{bdlem}), consider the  functions:
$$g(\Theta)=\Theta \quad \mbox{  and }  \quad g_\epsilon(\Theta) = C(\epsilon +
\epsilon\Theta + \Theta^2 + \Theta^4), $$ 
where $C$ is a given constant
and $\epsilon >0$ is a small parameter.

Clearly, $g(0) < g_\epsilon (0)$ for every $\epsilon$. Take now:
\begin{equation}\label{epsi}
\epsilon < \min\{\frac{1}{16C^2}, \frac{1}{4C}, 1\}
\end{equation}
and let $\Theta_0 \in (4C\epsilon, \frac{1}{4C})$ with
$\Theta_0<1$. Then: $\max\{C\epsilon, C\epsilon^2\Theta_0, C\Theta_0^2,
C\Theta_0^4\} <\frac{\Theta_0}{4}$ and hence $g(\Theta_0) > g_\epsilon(\Theta_0)$.

Taking now $T_{00}$ so small that, in addition to other requirements
imposed in the course of the proof, $\epsilon = T^{1/p} + D$ satisfies
(\ref{epsi}), we obtain that for every $T\in [0, T_{00})$ the
quantity $\Theta(T)$ must stay below $\Theta_0$, in virtue of  continuity of the
function $T\mapsto\Theta(T)$ and $\Theta(0) = 0$. This ends the proof
of (\ref{bdlem}) and of Lemma \ref{lemik}.
\endproof

\section{A proof of Theorem \ref{th1}}\label{mainproof}

We only outline the proof of Theorem \ref{th1}, which is standard, and
we point to its most important steps.
Let $\bar\xi$ be as in (\ref{barxi}). Recall that the system (\ref{v_eq}) can be
rewritten as:
\begin{equation}\label{sth}  
(\xi-\bar\xi)_{tt} - \mbox{div}
\left(D_Q\mathcal{Z}(\nabla\bar\xi,\nabla\bar\xi_t)\nabla(\xi-\bar\xi)_t\right)
= F[\xi,\bar\xi],
\end{equation}
where the right hand side $F[\xi,\bar\xi]$ is given in (\ref{Feq}).
We shall seek a solution $\xi$ as the fixed point of the operator:
$$\mathcal{T}(\tilde\xi - \bar\xi) = \xi - \bar\xi, \qquad \xi \mbox{
  is a solution to:}$$
\begin{equation}\label{sth2}
(\xi-\bar\xi)_{tt} - \mbox{div}
\left(D_Q\mathcal{Z}(\nabla\bar\xi,\nabla\bar\xi_t)\nabla(\xi-\bar\xi)_t\right)
= \tilde F[\tilde\xi,\bar\xi],
\end{equation}
in the Banach space:
\begin{equation}\label{space}
\begin{split}
E_{\Omega, T} = \Big\{u\in L_p(\Omega\times (0,T)); & ~ u(0,\cdot) = 0, 
~ u_t(0,\cdot) = 0,  ~u_{|\partial\Omega\times (0, T)} = 0, \\ &
\qquad\qquad \qquad\qquad ~ u_{tt}, \nabla^2u_t\in L_p(\Omega\times (0,T))\Big\},
\end{split}
\end{equation}
equipped with the norm:
$$\|u\|_{E_{\Omega, T}} = \Theta[u](T) = \|u_{tt}, \nabla^2
u_t\|_{L^p(\Omega\times (0,T)}.$$

\medskip

{\bf 1.} Following calculations as in the proof of Lemma \ref{lemik},
it results that:
$$\forall \tilde\xi - \bar\xi\in E_{\Omega, T} \qquad F[\tilde \xi,
\bar\xi]\in  L^p(\Omega\times (0,T).$$

\medskip

{\bf 2.} Integrating (\ref{sth2}) against $(\xi-\bar\xi)_t$ on
$\Omega\times (0,T)$ and using the estimate (\ref{positive}) in Lemma \ref{lemenergy}
below with $\zeta = (\xi-\bar\xi)_t$, we obtain:
\begin{equation*}
\begin{split}
\sup_{t\in (0,T)} \|(\xi-&\bar\xi)_t(t,\cdot)\|_{L_2(\Omega)}^2 + 
\|\nabla(\xi-\bar\xi)_t\|_{L_2(\Omega\times (0,T))}^2 \\ & \leq
C\int_0^T\int_\Omega F^2[\tilde\xi, \bar\xi]~\mbox{d}x ~\mbox{d}t 
+ C \|(\xi-\bar\xi)_t\|_{L_2(\Omega\times (0,T))}^2  \\ & \leq 
C\int_0^T\int_\Omega F^2[\tilde\xi, \bar\xi]~\mbox{d}x ~\mbox{d}t  
+ CT \sup_{t\in (0,T)} \|(\xi-\bar\xi)_t(t,\cdot)\|_{L_2(\Omega)}^2 ,
\end{split}
\end{equation*}
which implies the following energy estimate, for $T$ small:
\begin{equation}\label{energy}
\begin{split}
\sup_{t\in (0,T)} \|(\xi-\bar\xi)_t(t,\cdot)\|_{L_2(\Omega)}^2 + 
\|\nabla(\xi-\bar\xi)_t\|_{L_2(\Omega\times (0,T))}^2 \leq
C\int_0^T\int_\Omega F^2[\tilde\xi, \bar\xi]~\mbox{d}x ~\mbox{d}t.
\end{split}
\end{equation}
In virtue of (\ref{energy}), the Galerkin construction of the
approximants:
$$(\xi_N - \bar\xi)_t = \sum_{k=1}^N a_N^k(t) w_l(x),$$
where $\{w_l\}_{l=1}^\infty$ is an orthonormal base of $W_2^1(\Omega)$,
yields existence of a weak solution $\xi-\bar\xi = \lim_{N\to\infty}(\xi_N-\bar\xi)$
of the problem (\ref{sth2}), with: $(\xi-\bar\xi)_t\in L_\infty((0,T),
L_2(\Omega))$ and $\nabla (\xi-\bar\xi)_t\in L_2(\Omega\times (0,T)).$

\medskip

{\bf 3.} A modification of arguments in section \ref{secest} implies
that the weak solution $\xi$ is actually regular in the class
determined by (\ref{space}), i.e:
$$\xi-\bar\xi\in E_{\Omega, T}.$$
Moreover, for every small $\epsilon>0$:
$$ \mbox{if } \Theta[\tilde \xi - \bar\xi](T)\leq \epsilon \quad
\mbox{ then } \quad \Theta[\xi - \bar\xi](T)\leq \epsilon.$$ 

\medskip

{\bf 4.} In now suffices to show that the map $\mathcal{T}$ is a
contraction in some ball $\bar B_\epsilon\subset E_{\Omega, T}$. This
is done by applying methods of (\ref{secest}) to the system:
$$(\xi_1 - \xi_2)_{tt} - \mbox{div}
\left(D_Q\mathcal{Z}(\nabla\bar\xi,\nabla\bar\xi_t)\nabla(\xi_1-\xi_2)_t\right)
= F[\tilde\xi_1,\bar\xi] - F[\tilde\xi_2,\bar\xi],$$
where $\mathcal{T}(\tilde\xi_i - \bar\xi) = \xi_i - \bar\xi$. For
$\epsilon>0$ sufficiently small it follows that:
$$  \Theta[\xi_1 - \xi_2](T) \leq \frac{1}{2} \Theta[\tilde \xi_1 - \tilde\xi_2](T),$$
which completes the proof.
\endproof

\medskip

The key role above was played by the following estimate:
\begin{lemma}\label{lemenergy}
Let $T<T_0$ be sufficiently small and assume that $\mathcal{Z}$
satisfies (\ref{Zgood}). Then for every $\zeta\in
W_2^{2,1}(\Omega\times (0,T))$ such that $\zeta(0, \cdot)=0$ and
$\zeta_{|\partial\Omega\times (0,T)}=0$, there holds:
\begin{equation}\label{positive}
\begin{split}
\|\nabla \zeta\|^2_{L_2(\Omega\times (0,T))} \leq 4\gamma
\int_0^T\int_\Omega D_Q\mathcal{Z} (\nabla\bar\xi, \nabla\bar\xi_t)
\nabla\zeta : \nabla\zeta  ~\mathrm{d}x~\mathrm{d}t + 
C \|\zeta \|^2_{L_2(\Omega\times (0,T))}, 
\end{split}
\end{equation}
with constant $C$ independent of $\zeta$.
\end{lemma}
\begin{proof}
Consider a covering $\{B_k\}_{k=1}^N$ of $\Omega$ by a finite number
$N=N(r)$ of balls $B_k = B(X_k, r)$ with centers $X_k\in \Omega$ and
radius $r>0$. This family of coverings (parametrized by $r$) 
should be such that their covering numbers are uniform in $r$.
Let $\{\pi_k\}_{k=1}^N$ be a partition of unity subject to  
$\{B_k\}$.

For a fixed $t\in (0,T)$, with a slight abuse of notation, we shall
still write $\zeta = \zeta(t, \cdot) \in W^1_2(\Omega)$. By
(\ref{Zgood}) it follows that: 
\begin{equation}\label{f1}
\begin{split}
& \int_\Omega \langle D_Q\mathcal{Z}(\nabla\xi_0, \nabla
\xi_1)\nabla\zeta:\nabla\zeta\rangle \\ & = \sum_{k=1}^N \int_{B_k} 
\langle D_Q\mathcal{Z}(\nabla\xi_0(X_k), \nabla\xi_1(X_k))
\nabla(\pi_k^{1/2}\zeta):\nabla(\pi_k^{1/2}\zeta)\rangle~\mbox{d}x
+ \sum \int_{B_k} E_k[\xi, \bar\xi]  \\ & 
\geq \frac{1}{\gamma} \sum_{k=1}^N \|\nabla
(\pi_k^{1/2}\zeta)\|_{L_2(B_k)}^2 + \sum \int_{B_k} E_k[\xi, \bar\xi], 
\end{split}
\end{equation}
where we accumulated the error terms in:
\begin{equation*}
\begin{split}
E_k[\xi,\bar\xi] = & \langle D_Q\mathcal{Z}(\nabla\xi_0, \nabla
\xi_1)\pi_k^{1/2}\nabla\zeta:\pi_k^{1/2}\nabla\zeta \rangle \\ & - 
\langle D_Q\mathcal{Z}(\nabla\xi_0, \nabla
\xi_1)\nabla(\pi_k^{1/2}\zeta):\nabla(\pi_k^{1/2}\zeta) \rangle \\ & 
+ \langle [D_Q\mathcal{Z}(\nabla\xi_0, \nabla
\xi_1) - D_Q\mathcal{Z}(\nabla\xi_0(X_k), \nabla\xi_1(X_k))]
\nabla(\pi_k^{1/2}\zeta):\nabla(\pi_k^{1/2}\zeta) \rangle.
\end{split}
\end{equation*}
Hence:
\begin{equation*}
\begin{split}
|\int_{B_k}E_k&[\xi,\bar\xi]~\mbox{d}x| \leq   |\int_{B_k}\langle D_Q\mathcal{Z}(\nabla\xi_0, \nabla
\xi_1)\pi_k^{1/2}\nabla\zeta:(\zeta\otimes \nabla\pi_k^{1/2}) \rangle|
\\ & +
|\int_{B_k}\langle D_Q\mathcal{Z}(\nabla\xi_0, \nabla \xi_1)
(\zeta\otimes \nabla\pi_k^{1/2}) : \nabla(\pi_k^{1/2}\zeta) \rangle| +
C r \|\nabla(\pi_k^{1/2}\zeta) \|_{L_2(B_k)}^2 \\ &
\leq C_r \|\zeta \|_{W^1_2(B_k)} \|\zeta\|_{L_2(B_k)} + Cr \|\nabla(\pi_k^{1/2}\zeta) \|_{L_2(B_k)}^2,   
\end{split}
\end{equation*}
where $C$ is a universal constant depending only on the initial data
and $\mathcal{Z}$, while the constant $C_r$ depends on the covering $\{B_r\}$.
Taking $r$ small, so that $Cr<1/(2\gamma)$, by (\ref{f1}) we now arrive at:
\begin{equation*}
\begin{split}
 \int_\Omega \langle D_Q\mathcal{Z}(\nabla\xi_0, &\nabla
\xi_1)\nabla\zeta:\nabla\zeta\rangle \\ & \geq \frac{1}{2\gamma} \sum_{k=1}^N \|\nabla
(\pi_k^{1/2}\zeta)\|_{L_2(B_k)}^2 - C_r \|\zeta
\|_{W^1_2(\Omega)} \|\zeta\|_{L_2(\Omega)} \\ & 
\geq \frac{1}{2\gamma} \sum_{k=1}^N \|\pi_k^{1/2}\nabla\zeta\|_{L_2(B_k)}^2 - C_r \|\zeta
\|_{W^1_2(\Omega)} \|\zeta\|_{L_2(\Omega)} \\ &
\geq \frac{1}{2\gamma} \|\nabla\zeta\|_{L_2(\Omega)}^2 - C_r \left(\epsilon \|\nabla\zeta
\|^2_{L_2(\Omega)} + \frac{1}{\epsilon}\|\nabla\zeta\|^2_{L_2(\Omega)} \right),
\end{split}
\end{equation*}
by Young's inequality. With $\epsilon$ sufficiently small, it yields:
\begin{equation}\label{f2}
\|\nabla\zeta\|_{L_2(\Omega)}^2 \leq 3\gamma 
 \int_\Omega \langle D_Q\mathcal{Z}(\nabla\xi_0, \nabla
\xi_1)\nabla\zeta:\nabla\zeta\rangle + C \|\zeta\|^2_{L_2(\Omega)}.
\end{equation}
Integrating in $t$, we eventually arrive at:
\begin{equation*}
\begin{split}
\|\nabla \zeta\|_{L_2(\Omega\times (0,T))}^2 & \leq 3\gamma 
\int_0^T \int_\Omega \langle D_Q\mathcal{Z}(\nabla\xi_0, \nabla
\xi_1)\nabla \zeta:\nabla \zeta\rangle ~\mbox{d}x~\mbox{d}t
+ C \|\zeta\|^2_{L_2(\Omega\times (0, T))} \\ &
\leq 3\gamma 
\int_0^T \int_\Omega \langle D_Q\mathcal{Z}(\nabla\bar\xi, \nabla
\bar\xi_t)\nabla \zeta :\nabla \zeta\rangle
~\mbox{d}x~\mbox{d}t
\\ & \qquad \qquad \qquad+ CT \| \nabla \zeta\|^2_{L_2(\Omega\times (0, T))} 
+ C \| \zeta\|^2_{L_2(\Omega\times (0, T))},
\end{split}
\end{equation*}
which for $T$ small enough implies (\ref{positive}).
\end{proof}

\section{A proof of Lemma \ref{Zexamples}}\label{seclast}

{\bf 1.} To prove (i), note that $D_Q\mathcal{Z}_0''(F_0, Q_0) Q =
2F_0\mbox{sym}(F_0^TQ)$ so that:
$$\forall Q\in\mathbb{R}^{n\times n} 
\qquad \langle D_Q\mathcal{Z}_0''(F_0, Q_0) Q :Q\rangle= 2 \langle
\mbox{sym}(F_0^TQ_0): F_0^TQ\rangle = |\mbox{sym}(F_0^TQ)|^2.$$
Take $\zeta\in W_2^1(\Omega,\mathbb{R}^n)$ with trace $0$ on the boundary $\partial \Omega$.
We have:
\begin{equation*}
\begin{split}
\int_\Omega|\nabla\zeta|^2 & \leq |F_0^{-1,T}|^2 \int_\Omega|\nabla(F_0^T\zeta)|^2
\leq 2 |F_0^{-1,T}|^2 \int_\Omega|\mbox{sym}\nabla(F_0^T\zeta)|^2
\\ & = |F_0^{-1,T}|^2 \int_\Omega \langle D_Q\mathcal{Z}_0''(F_0, Q_0)
\nabla\zeta: \nabla\zeta\rangle,
\end{split}
\end{equation*}
where we applied Korn's inequality to the map $x\mapsto F_0^T\zeta(x)$.

\medskip

{\bf 2.} To prove (ii), observe that $D_Q\mathcal{Z}_0'(F_0, Q_0) Q =
2(\det F_0) \mbox{sym}(QF_0^{-1})F_0^{-1,T}$ so that:
$$\langle D_Q\mathcal{Z}_0'(F_0, Q_0) Q :Q\rangle= 2 (\det F_0) |\mbox{sym}(QF_0^{-1})|^2.$$
Then, for any test function $\zeta$ as above, we have:
\begin{equation}\label{c2}
\begin{split}
\int_\Omega|\nabla\zeta|^2 & \leq |F_0|^2 \int_\Omega|(\nabla\zeta) F_0^{-1}|^2~\mbox{d}x
= |F_0|^2 \int_{F_0 \Omega}|\nabla (\zeta\circ (F_0^{-1}y))|^2(\det
F_0^{-1})~\mbox{d}y \\ & \leq  2 |F_0|^2 (\det F_0^{-1})
\int_{F_0 \Omega} |\mbox{sym}((\nabla \zeta)\circ F_0^{-1})F_0^{-1})|^2~\mbox{d}y
\\ & = 2|F_0|^2 (\det F_0^{-1})(\det F_0)
\int_\Omega |\mbox{sym}((\nabla \zeta)F_0^{-1})|^2~\mbox{d}x  \\ & = 
|F_0|^2 (\det F_0)^{-1}\int_\Omega \langle D_Q\mathcal{Z}_0'(F_0, Q_0)
\nabla\zeta: \nabla\zeta\rangle,
\end{split}
\end{equation}
where we applied Korn's inequality to the map $y\mapsto
\zeta(F_0^{-1}y)$ on the open domain $F_0\Omega$.

\medskip

{\bf 3.} To prove (iii) -- (v), observe that:
\begin{equation*}
\begin{split}
\langle D_Q\mathcal{Z}_m(F_0, Q_0) Q :Q\rangle & = 
\Big\langle \sum_{j=0}^{2m}  (\mbox{sym}(Q_0F_0^{-1})^j 
\mbox{sym} (QF_0^{-1}) (\mbox{sym}(Q_0F_0^{-1})^{2m-j} : QF_0^{-1}
\Big\rangle \\ & = \Big\langle \sum_{j=0}^{2m} A^j B A^{2m-j} 
: QF_0^{-1}\Big\rangle,
\end{split}
\end{equation*}
where we denoted: 
$$A = \mbox{sym}(Q_0F_0^{-1}), \qquad B = \mbox{sym}(QF_0^{-1}).$$
Since the matrix $\sum_{j=0}^{2m} A^j B A^{2m-j} $ is symmetric, it
follows that:
$$\langle D_Q\mathcal{Z}_m(F_0, Q_0) Q :Q\rangle= 
\Big\langle \sum_{j=0}^{2m} A^j B A^{2m-j} : B\Big\rangle$$

Let $\zeta$ be a test function as in Lemma \ref{lem1}. By calculations
similar to (\ref{c2}) we get:
$$\int_\Omega |\nabla\zeta|^2 \leq \frac{1}{2}|F_0|^2 
\int_\Omega \langle D_Q\mathcal{Z}_0(F_0, Q_0)
\nabla\zeta: \nabla\zeta\rangle, $$
proving (iii). To prove (iv), we compute:
\begin{equation*}
\begin{split}
\langle D_Q\mathcal{Z}_1&(F_0, Q_0) Q :Q\rangle = 
\langle A^2B : B\rangle + \langle ABA : B\rangle + \langle BA^2 : B\rangle 
\\ & = \langle AB : AB\rangle \langle BA : AB\rangle  + |AB|^2
= 2 \langle \mbox{sym}(AB) : AB\rangle + |AB|^2 \\ &
= 2|\mbox{sym}(AB)|^2 + |AB|^2  \geq |AB|^2.
\end{split}
\end{equation*}
Therefore, by calculations similar to (\ref{c2}):
 \begin{equation}\label{c4}
\begin{split}
\int_\Omega|\nabla\zeta|^2 & \leq 2|F_0|^2 
\int_\Omega |\mbox{sym}((\nabla \zeta)F_0^{-1})|^2  
\leq 2|F_0|^2 |A^{-1}|^2 \int_\Omega |AB|^2 \\ & \leq 2|F_0|^2
|\mbox{sym}(Q_0F_0^{-1})^{-1}|^2 
\int_\Omega \langle D_Q\mathcal{Z}_1(F_0, Q_0)
\nabla\zeta: \nabla\zeta\rangle.
\end{split}
\end{equation}
Finally, in order to prove (v) we derive:
\begin{equation*}
\begin{split}
\langle D_Q\mathcal{Z}_2(F_0, Q_0) Q :Q\rangle & = 
\langle A^4B  + A^3BA + A^2BA^2 + ABA^3 + BA^4: B\rangle
\\ & = |A^2B + ABA|^2 + |A^2B|^2\geq |A^2B|^2,
\end{split}
\end{equation*}
which, in the same manner as in (\ref{c4}) yields:
 \begin{equation*}
\begin{split}
\int_\Omega|\nabla\zeta|^2 & \leq 2|F_0|^2 
|A^{-2}|^2 \int_\Omega |A^2B|^2 \\ & \leq 2|F_0|^2
|\mbox{sym}(Q_0F_0^{-1})^{-1}|^4
\int_\Omega \langle D_Q\mathcal{Z}_2(F_0, Q_0)
\nabla\zeta: \nabla\zeta\rangle.
\end{split}
\end{equation*}
The proof of Lemma \ref{Zexamples} is done.

\end{document}